\documentclass{amsart}

\usepackage{amssymb} % \mathbb
\usepackage{amsthm}
\usepackage{amsmath}
\usepackage{latexsym}
\usepackage{mathrsfs}
\usepackage{enumerate}
\usepackage[dvips]{graphicx}
\usepackage{color}
\newtheorem{theorem}{Theorem}[section]

\newtheorem{example}[theorem]{Example}

\definecolor{mgre}{cmyk}{0.92,0.00,0.59,0.25}

\def\bR{{\mathbb{R}}}

\def\bfE{{\mathbf{E}}}

\def\bfM{{\mathbf{M}}}

\def\bfP{{\mathbf{P}}}

\def\cD{{\mathcal{D}}}
\def\cE{{\mathcal{E}}}
\def\cF{{\mathscr{F}}}

\def\cH{{\mathcal{H}}}

\def\cK{{\mathcal{K}}}
\def\cL{{\mathcal{L}}}

\def\cS{{\mathcal{S}}}

\def\sB{{\mathscr{B}}}

\def\1{{\mathbf{1}}}

\newcommand{\capa}{\mathrm{Cap}}

\def\<{\langle}
\def\>{\rangle}

\newtheorem{thm}{Theorem}[section]
\newtheorem{lem}{Lemma}[section]
\newtheorem{cor}{Corollary}[section]
\newtheorem{prop}{Proposition}[section]
\newtheorem{rem}{Remark}[section]

\newtheorem{defn}{Definition}[section]
\newtheorem{pro}{Problem}[section]
\newcommand{\supp}{\textrm{supp}}

\newcommand{\bd}{\begin{defi}}
\newcommand{\ed}{\end{defi}}

\newcommand{\bpro}{\begin{pro}}
\newcommand{\epro}{\end{pro}}

\newcommand{\bec}{\begin{cases}}
\newcommand{\eec}{\end{cases}}

\newcommand{\bpr}{\begin{prob}}
\newcommand{\epr}{\end{prob}}

\newcommand{\bt}{\begin{thm}}
\newcommand{\et}{\end{thm}}

\newcommand{\ba}{\begin{ass}}
\newcommand{\ea}{\end{ass}}

\newcommand{\br}{\begin{rem}}
\newcommand{\er}{\end{rem}}

\newcommand{\bpm}{\begin{pmatrix}}
\newcommand{\epm}{\end{pmatrix}}

\newcommand{\be}{\begin{ex}}
\newcommand{\ee}{\end{ex}}

\newcommand{\bp}{\begin{prop}}
\newcommand{\ep}{\end{prop}}

\newcommand{\bl}{\begin{lem}}
\newcommand{\el}{\end{lem}}

\newcommand{\bc}{\begin{cor}}
\newcommand{\ec}{\end{cor}}

\newcommand{\bq}{\begin{que}}
\newcommand{\eq}{\end{que}}

\newcommand{\beqn}{\begin{eqnarray*}}
\newcommand{\eeqn}{\end{eqnarray*}}

\newcommand{\beqnn}{\begin{eqnarray}}
\newcommand{\eeqnn}{\end{eqnarray}}

\newcommand{\bequ}{\begin{equation}}
\newcommand{\eequ}{\end{equation}}

\newcommand{\benu}{\begin{enumerate}}
\newcommand{\eenu}{\end{enumerate}}

\newcommand{\barr}{\begin{array}{rcl}}
\newcommand{\ear}{\end{array}}

\newcommand{\la}{\label}

\newcommand{\eps}{\epsilon}

\newcommand{\Om}{\Omega}

\usepackage{amsmath}	% required for `\cases' (yatex added)

\begin{document}

\title [the maximum principle of Schr\"{o}dinger operators]{The bottom of the spectrum of time-changed processes and the maximum principle of Schr\"{o}dinger operators}

\author{Masayoshi Takeda}
\address{Mathematical Institute,
Tohoku University, Aoba, Sendai, 980-8578, Japan}
\email{takeda@math.tohoku.ac.jp}
\thanks{The author was supported in part by Grant-in-Aid for Scientific
Research (No.26247008(A)) and Grant-in-Aid for Challenging Exploratory Research (No.25610018), Japan Society for the Promotion of Science.}

%    General info
%\date{\today}


\begin{abstract}
We give a necessary and sufficient condition for 
the maximum principle of Schr\"{o}dinger operators in terms of 
the bottom of the spectrum of time-changed processes. 
As a corollary, we obtain a sufficient condition for the Liouville property of Schr\"{o}dinger operators. 
\end{abstract}

\maketitle

\section{Introduction}
In \cite{T3}, we define the {subcriticality}, {criticality} and {superciriticality} for   
Schr\"{o}dinger forms and characterize these properties in terms of the bottom of the spectrum of time
changed processes. In the process, we prove the existence of a harmonic function (or ground state) of the 
 Schr\"odinger form and study its properties. In particular, we show that it has a bounded, positive,  
continuous version which is invariant with respect to its Schr\"odinger semigroup. 
In this paper, we will show, as an application of this fact, the maximum 
principle and Liouville property of Schr\"odinger operators.

Let $X$ be a locally compact separable metric space and $m$ a positive 
Radon measure on $X$ with full topological support. Denote by $X_\triangle :=X\cup \{\triangle \}$ the one-point 
compactification of $X$.
Let $\bfM = (\bfP_x, X_t,\zeta )$ be an $m$-symmetric Hunt process with lifetime $\zeta=\inf\{t>0\mid X_t=\triangle \}$. 
We assume that 
$\bfM$ is irreducible and strong Feller.
% and, in addition, that 
%${\bf M}$ generates a {\it regular} Dirichlet form $(\cE, \cD(\cE))$ on $L^2(X;m)$.  
Let $\mu=\mu^+-\mu^-$ be a signed Radon smooth measure such that the positive (resp. negative) part $\mu^+$ (resp. 
$\mu^-$) belongs to 
the local Kato class (resp. the Kato class). 
%Throughout this paper, we make the convention that any function $f$ defined on $X$ is extended 
%to $X_\triangle $ by setting $f(\triangle )=0$.
%In \cite{T3}, we define the criticality or subcriticality for  $\cE^\mu$ 
%through $h$-transform; 
We denote by $A^{\mu^+}_t$ (resp. $A^{\mu^-}_t$) the positive continuous additive functional 
in the Revuz correspondence to $\mu^+$ (resp. $\mu^+$). Put $A^{\mu}_t=A^{\mu^+}_t-A^{\mu^-}_t$
and define the Feynman-Kac semigroup $\{p^\mu_t\}_{t\geq 0}$ by
$$
p^\mu_tf(x)=\bfE_x\left(e^{-A^\mu_t}f(X_t)\right).
$$

%We denote by $\cD_{\loc}(\cE)$ 
%the set of functions locally in $\cD(\cE)$ and introduce a function space by
%\begin{equation*}
%\cH^+(\mu) = \{ h \in \cD_{\loc}(\cE) \cap C(X) : h>0,\  p_t^{\mu} h \le h\}.
%\end{equation*}
%Namely, $\cH^+(\mu)$ is the space of $p^\mu_t$-excessive functions. 

We denote by ${\bf M}^{\mu^+}=(P^{\mu^+}_x,X_t,\zeta )$ the subprocess of ${\bf M}$ by the multiplicative 
functional $\exp({-A^{\mu^+}_t})$ and by $(\cE^{\mu^+},\cD(\cE^{\mu^+}))$ the Dirichlet
 form generated by ${\bf M}^{\mu^+}$. Suppose that 
 the negative part $\mu^-$ is non-trivial and Green-tight 
with respect to ${\bf M}^{\mu^+}$ (Definition \ref{def-Kato} (2)). 
We then define
  $\lambda (\mu)$ by
\begin{equation}\la{i-ground}
\lambda (\mu) := \inf \left\{ \cE(u,u) + \int_X u^2d\mu^+\mid u \in
		       \cD(\cE),\ \int_X u^2 d\mu^- = 1\right\}.
\end{equation}
$\lambda (\mu)$ is regarded as the bottom of the spectrum of the time-changed process 
of ${\bf M}^{\mu^+}$ by the continuous additive functional $A^{\mu^-}_t$. We show in \cite[Theorem 2.1]{T2} that the minimizer of \eqref{i-ground} exists in the extended 
Dirichlet space $\cD_e(\cE^{\mu^+})$ and it can be taken to be strictly positive on $X$. 
 The objective of this paper is to prove 
 the maximum principle of Schr\"odinger forms by using the existence of 
 the minimizer of (\ref{i-ground}).
More precisely, let 
\bequ\la{ba}
\cH^{ba}(\mu)= \{ h\in\sB(X)\mid h\ \text{ is bounded above},\  p_t^{\mu} h \geq  h\},
\eequ
where $\sB(X)$ is the set of Borel functions on $X$. We here define the maximum principle by

\bigskip
\noindent
\ \ \ \ ({\bf MP})\ \  If $h\in \cH^{ba}(\mu)$, then $h(x)\leq 0$ for all $x\in X$.

\bigskip
\noindent
We will prove in Theorem \ref{MP} that under Assumption

\bigskip
\noindent
\ \ \ \ {\rm({\bf A})}\ \  $\bfE_x\left(e^{-A^{\mu^+}_\infty };\zeta =\infty \right)=0$,

\bigskip
\noindent
{\rm({\bf MP})} is equivalent to $\lambda (\mu )>1$. 
For the proof of this, it is crucial that if $\lambda( \mu ) = 1$, then  
the minimizer $h$ in (\ref{i-ground}) has a bounded continuous version with $p_t^\mu$-invariance, i.e., 
$\bfE_x(\exp(-A^\mu_t)h(X_t))=h(x)$ (\cite[Lemma 5.16, Corollary 5.17]{T3}).

\bigskip
Let us introduce the space $\cH^b(\mu)$ of bounded $p_t^\mu$-invariant functions:
$$
\ \cH^b(\mu)= \{ h\in\sB_b(X) \mid  p_t^{\mu} h = h\}.
$$
We here define the {\it Liouville property} by

\bigskip
\noindent
\ \ \ \ ({\bf L})\ \ If $h\in {\cH}^{b}(\mu)$, then $h(x)=0$ for all $x\in X$.

\bigskip
\noindent
We will show in Corollary \ref{Cor} that under Assumption ({\bf A}), 
$\lambda (\mu )>1$ implies ({\bf L}).

\bigskip
We remark that Theorem \ref{MP} and Corollary \ref{Cor} can be applied to non-local Dirichlet forms.  
In a remaining part of introduction, we treat these two properties  
for strongly local Dirichlet forms, which are regarded as an extension of symmetric 
elliptic operators of second order. In Berestycki-Nirenberg-Varadhan \cite{BNV}, 
they define a {maximum principle} for a uniformly elliptic operator of second order,
$L=M+c=a_{i,j}\partial _i\partial _j+b_i\partial _i+c$,  
on a general bounded domain $D$ of ${\mathbb R}^d$. Let $u_0$ be the solution to the equation $Mu=-1$ 
vanishing on $\partial D$ in a suitable sense: define $\cS$ by the set of sequences $\{x_n\}_{n=1}^\infty \subset D$ such that
$x_n$ converges to a point of the boundary $\partial D$ and $u_0(x_n)$ converges to 0. 
They say that the {\it refined maximum principle} holds for $L$, 
if a function $h$ bounded above satisfies $Lh \geq 0$ on $D$ and $\limsup_{n\to\infty }h(x_n)\leq 0$ for any $\{x_n\}_{n=1}^\infty
\in\cS$, then $h\leq 0$ on $D$, and prove that $L$ satisfies the refined maximum principle if and only if
the principal eigenvalue $\lambda _0$ of $-L$ is positive. 

Note that $u_0$ equals $\bfE_x(\tau_D)$, where $\bfP_x$ is
 the diffusion process with generator
$M$ and $\tau_D$ is the first exit time from $D$. We see that if $D$ is bounded (more generally, 
Green-bounded, i.e., $\sup_{x\in D}\bfE_x(\tau_D)<\infty $), then $\cS$ is identical to
 the set of sequences $\{x_n\}$ such that $x_n\to\partial D$ and $\bfE_{x_n}\left(\exp(-\tau_D)\right)\to 1$
 as $n\to\infty $ (Lemma \ref{S=S}, Remark \ref{S=S1}). 
 Considering this fact, we define 
\begin{equation}\la{cS}
\cS=\left\{\{x_n\}_{n=1}^\infty\subset X\mid 
\begin{split}
&\text{$x_n\to \triangle $ and $\bfE_{x_n}(e^{-\zeta })\to 1$ as $n\to\infty $}\\
\end{split}
\right\}.
\end{equation}
Assume that $(\cE,\cD(\cE))$ is {strongly local} and set
\begin{equation*}
\widetilde{\cH}^{ba}(\mu)=\left\{h \mid
\begin{split}
&\text{$h\in\cD_{loc}(\cE)\cap C(X)$ is bounded above, $\cE^\mu(h,\varphi )\leq 0$ for}\\
&\text{$\forall \varphi \in \cD(\cE)\cap C^+_0(X)$ and $\limsup_{n\to\infty }h(x_n)\leq 0$ 
for $\forall \{x_n\} \in \cS$}
\end{split}
\right\},
\end{equation*}
where $C^+_0(X)$ is the set of non-negative continuous functions with compact support. 
Following \cite{BNV}, we here define the {\it refined maximum principle} by

\bigskip
\noindent
\ \ \ \ (${\bf{RMP}}$)\ \ If $h\in \widetilde{\cH}^{ba}(\mu)$, then $h(x)\leq 0$ for all $x\in X$.

\bigskip
\noindent
We will show that $\widetilde{\cH}^{ba}(\mu)\subset {\cH}^{ba}(\mu)$ (Lemma \ref{include}), and thus see, as a corollary of Theorem of \ref{MP}, that  
$\lambda (\mu)>1$ implies ({\bf RMP}) (Theorem \ref{RMP}).
%If the minimizer $\phi $ in (\ref{i-ground}) satisfies $\limsup_{n\to\infty }\phi (x_n)\leq 0$ for any $\{x_n\} \in \cS$, then the converse, ({\bf RMP}) implies $\lambda (\mu)>1$, also holds. Since the decay of $\phi $ equals that of the Green function $G^{\mu^+}(o,\cdot )$
% in some regular cases, this assumption is fulfilled if for a fixed point $o\in X$ $\lim_{x\to\triangle }G^{\mu^+}(0,x)=0$.
We would like to emphasize that if $D$ is bounded and $L$ is symmetric, the principal eigenvalue $\lambda _0$ 
of $-L$ is positive 
if and only if $\lambda (\mu)>1$.
However, $\lambda(\mu)>1$ does not always imply $\lambda_0>0$ for a unbounded domain $D$, while
$\lambda _0>0$ implies $\lambda (\mu )>1$ in general (Lemma \ref{first}).

\medskip
When $(\cE,\cD(\cE))$ is {strongly local}, we set 
\begin{equation*}
\widetilde{\cH}^b(\mu)=\{h\in\cD_{loc}(\cE)\cap C_b(X)\mid 
\cE^\mu(h,\varphi )=0, \forall \varphi \in \cD(\cE)\cap C_0(X)\}
\end{equation*}
and define the property $(\widetilde{\bf L})$ by

\bigskip
\noindent
\ \ \ \ $(\widetilde{\bf L})$\ \ If $h\in\widetilde{\cH}^b(\mu)$, then $h(x)=0$ for all $x\in X$.

\bigskip
\noindent
We then see that if ${\bfM}$ is conservative and ({\bf A}) is fulfilled, 
that is, $\bfE_x(\exp({-A^{\mu^+}_\infty }))=0$,
 then $\widetilde{\cH}^{b}(\mu)\subset {\cH}^{b}(\mu)$, and consequently 
 $\lambda (\mu )>1$ implies ($\widetilde{\bf L}$) 
 (Corollary \ref{Cor1}). Grigor'yan and Hansen \cite{GH} calls
 a measure $\mu^+$ {\it big} if it satisfies ({\bf A}), and they 
prove that for the transient Brownian motion 
${\bfM}=(\bfP_x,B_t)$ on $\mathbb{R}^d$, if $\mu^-\equiv 0$ and  $\mu^+$ is big, 
then ($\widetilde{\bf L}$) holds. Corollary \ref{Cor} tells us that 
if $\mu^-$ is small with respect to $\mu^+$ in the sense that $\lambda (\mu )>1$
then $(\widetilde{\bf L})$ still holds.

Pinsky \cite{RP} treat absolutely continuous potentials $d\mu=V^+dx-V^-dx$ and prove in \cite[Theorem 1.1]{RP}
that if $
\sup_{x\in{\mathbb R}^d}\bfE_x\left(\exp({\int_0^\infty V^-(B_t)dt})\right)<\infty,
$
the Liouville property $(\widetilde{\bf L})$ is equivalent to 
$$
\int_0^\infty V^+(B_t)dt=\infty,\ \bfP_x\text{-a.e.}\ (\Longleftrightarrow \ ({\bf A})).
$$ 
We will give an example of potential $\mu$ that even if 
$\sup_{x\in{\mathbb R}^d}\bfE_x(\exp(A^{\mu^-}_\infty))=\infty$ and $\bfE_x(\exp(-A^{\mu^+}_\infty))=0$, 
 $(\widetilde{\bf L})$ holds (Example \ref{e-2}).

\section{Schr\"odinger forms}
Let $X$ be a locally compact separable metric space and $m$ a positive
Radon measure on $X$ with full topological support.
Let $(\cE, \cD(\cE))$ be a regular Dirichlet form on
$L^2(X;m)$. 
We denote by $u\in \cD_{loc}(\cE)$ if for any relatively compact open set $D$ there exists a function $v\in\cD(\cE)$
such that $u=v$ $m$-a.e. on $D$. 
%Let $\cL$ be the self-adjoint operator associated with $\cE$. 
We denote by $\cD_e(\cE)$ the family of $m$-measurable functions $u$ on
$X$ such that $|u| < \infty$ $m$-a.e. and there exists an $\cE$-Cauchy
sequence $\{ u_{n}\}$ of functions in $\cD(\cE)$ such that
$\lim_{n \to \infty} u_{n} = u$ $m$-a.e.
We call $\cD_e(\cE)$ the
\emph{extended Dirichlet space} of $(\cE, \cD(\cE))$.

 Let $\bfM = (\Om, \cF, \{\cF_{t}\}_{t \ge 0}, \{\bfP_x\}_{x \in X}, \{X_t\}_{t \ge 0},
\zeta)$ be the symmetric Hunt process generated by $(\cE, \cD(\cE))$,
where $\{ \cF_t\}_{t \ge 0}$ is the augmented filtration and $\zeta $ is the lifetime of ${\bf M}$.
 Denote by $\{p_t\}_{t\geq 0}$ and  $\{G_\alpha\}_{\alpha \geq 0} $ the semigroup and resolvent 
of ${\bf M}$:
$$
p_tf(x)=\bfE_x(f(X_t)),\ \ \ \ \ G_\alpha f(x)=\int_0^\infty e^{-\alpha t} p_tf(x)dt.
$$

We assume that ${\bf M}$ satisfies next two conditions:  

\smallskip
%\noindent
 \begin{description}
	\item{\textbf{Irreducibility (I)}}. \ If a Borel set $A$ is $p_t$-invariant, i.e., 
      $p_t(1_Af)(x)=1_Ap_tf(x)$\ $m$-a.e. for any $f\in
      L^2(X;m)\cap{\mathscr B}_b(X)$  
      and $t>0$, then $A$ satisfies either $m(A)=0$ or $m(X\setminus A)=0$. 
      Here ${\mathscr B}_b(X)$ is the space of bounded Borel functions on
      $X$.
\item{\textbf{Strong Feller Property (SF)}}. \  For each $t$,
      $p_t ({\mathscr  B}_b(X))\subset C_b(X)$, 
      where $C_b(X)$ is the space of bounded continuous functions on $X$.
\end{description}
\smallskip

We remark that \textbf{(SF)} implies \textbf{(AC)}. 
\begin{description}
	\item{\textbf{Absolute Continuity Condition (AC)}}. \ The transition probability of $\bfM$ is absolutely continuous with respect to $m$, $p(t,x,dy)=p(t,x,y)m(dy)$ for each $t>0$ and $x\in X$.
 \end{description}

%\begin{rem}
% We see from the assumption ({\bfseries SF}) that the semigroup
% $\{P_{t}\}_{t>0}$ admits an integral kernel $\{p(t,x,y)\}_{t>0}$ with
% respect to the measure $m$.
%\end{rem}
Under {\textbf{(AC)}}, a non-negative, jointly measurable $\alpha$-resolvent kernel $G_{\alpha}(x,y)$ exists: 
$$
G_\alpha f(x)=\int_XG_\alpha (x,y)f(y)m(dy),\ x\in X,\ f\in{\mathscr B}_b(X).
$$
Moreover, $G_\alpha (x,y)$ is $\alpha $-excessive in $x$ and in $y$ (\cite[Lemma 4.2.4]{FOT}).
We  simply write $G(x,y)$ for $G_{0}(x,y)$. 
For a measure $\mu$, we define the $\alpha $-potential of $\mu$ by
\begin{equation*}
G_{\alpha} \mu(x) = \int_{X} G_{\alpha}(x,y) \mu(dy).
\end{equation*}
%If $\alpha = 0$, we define 
%\begin{equation*}
%G \mu(x) := G_0 \mu(x).
%\end{equation*}

\begin{defn}\label{def-transient}\rm
 \begin{enumerate}[(1)]
  \item A Dirichlet space $(\cE, \cD(\cE))$ on $L^2(X;m)$ is said to be \emph{transient} if there exists a strictly positive, bounded
	function $g$ in $L^1(X;m)$  such that for $u \in \cD(\cE)$
	\begin{equation*}
	 \int_X |u| g dm \le \sqrt{\cE(u,u)}. 
	\end{equation*}
%	(cf. \cite[p.40]{FOT}).
  \item  A Dirichlet space $(\cE, \cD(\cE))$ on $L^2(X;m)$ is said to be \emph{recurrent} if the constant function 1
belongs to $\cD_e(\cE)$ and $\cE(1,1)=0$. Namely, there exists a
	sequence $\{ u_n\} \subset \cD(\cE)$ such that $\lim_{n, m \to \infty} \cE(u_n - u_m , u_n - u_m) = 0$ and 
	%$0 \le u_n \le 1, 
$\lim_{n \to \infty} u_n  = 1 \ \textrm{$m$-a.e.}$	
 \end{enumerate}
\end{defn}
For other characterizations of transience and recurrence, see \cite[Theorem 1.6.2, Theorem 1.6.3]{FOT}.

We define the {\em (1-)capacity\/} $\capa$ associated with the Dirichlet
form $(\cE, \cD(\cE))$ as follows:
for an open set $O\subset X$,
\begin{equation*}
 \capa(O)=\inf\left\{\cE(u,u)+(u,u)_m\mid u\in \cD(\cE), u\geq 1\ m\text{-a.e. on
  }O\right\}
\end{equation*}
and for a Borel set $A\subset X$,
\begin{equation*}
 \capa(A)=\inf\{\capa(O)\mid O\text{ is open, }O\supset  A\}.
\end{equation*}
A statement depending on $x\in X$ is said to hold q.e. on $X$ if there
exists a set $N\subset X$ of zero capacity such that the statement is
true for every $x\in X\setminus N$.
``q.e.'' is an abbreviation of ``quasi-everywhere''.
A real valued function $u$ defined q.e. on $X$
is said to be {\em quasi-continuous\/} if for any
$\epsilon>0$ there exists 
an open set $G\subset X$ such that $\capa(G)<\epsilon$ and
$u|_{X\setminus G}$ is finite and  continuous.
Here, $u|_{X\setminus G}$ denotes the restriction of $u$ to $X\setminus
G$.
Each function $u$ in $\cD_e(\cE)$ admits a
quasi-continuous version $\tilde{u}$, that is,
$u=\tilde{u}$ $m$-a.e. % and $\tilde{u}$ is quasi continuous.
In the sequel,
we always assume that every function $u\in \cD_e(\cE)$ is
represented by its quasi-continuous version.

We call a positive Borel measure $\mu$ on $X$ {\it smooth} if it satisfies the following conditions:

\medskip
{(S1)} \ $\mu$ charges no set of zero capacity,

\medskip
(S2) \ there exists an increasing sequence $\{F_n\}$ of closed sets that
\bequ\la{s1}
\mu(F_n)<\infty,
\eequ
\bequ\la{s2}
\lim_{n\to\infty }{\rm Cap}(K\setminus F_n)=0\quad \text{for any compact set $K$}.
\eequ

\medskip
\noindent
We denote by $S$ the set of smooth measures.
%it satisfies
% the following conditions:
%\smallskip 
%(i)\ \ for all $A\in {\mathscr B}(X)$, $\capa(A)=0$ implies $\mu(A)=0$.
%\smallskip
%(ii) there exists an increasing sequence $\{F_{n}\}$ of  closed sets
%	such that 
%	$\mu(F_{n})<\infty$ for all $n$ and
%	  $\lim_{n\rightarrow \infty}\capa(K\setminus F_{n})=0$
%	  for any compact set $K$.
%\smallskip
%We denote by $\calS$ is the family of all smooth measures.

A stochastic process $\{A_{t}\}_{t\geq 0}$ is said to be an {\em additive
functional\/} %\label{defn:addtive functional}
(AF in abbreviation)
if the following conditions hold:
\begin{enumerate}[(i)]
\item  $A_{t}(\cdot)$ is ${\cF}_{t}$-measurable for all $t\geq 0$.

\item  There exists a set $\Lambda\in
       {\cF}_{\infty}=\sigma\left(\cup_{t\geq 0}{\cF}_{t}\right)$ such
       that 
       $P_{x}(\Lambda)=1$, for q.e. $x\in  X$,
       $\theta_{t}\Lambda\subset\Lambda$ for all $t>0$,
       and for each $\omega\in \Lambda$, $A_{\cdot}(\omega)$ is a
       function satisfying:
       $A_{0}=0$, $A_{t}(\omega)<\infty$ for $t<\zeta(\omega)$,
       $A_{t}(\omega)=A_{\zeta}(\omega)$ for $t\geq \zeta $,
       and $A_{t+s}(\omega)=A_{t}(\omega)+A_{s}(\theta_{t}\omega)$ for
       $s,t\geq 0$.
\end{enumerate}
\noindent
If an AF $\{A_{t}\}_{t\geq 0}$ is positive and continuous with respect
to $t$ for each $\omega\in \Lambda$, the AF is called a {\em positive
continuous additive functional\/} (PCAF in abbreviation). The set of all PCAF's is denoted by ${\bf A}^+_c$. 
The family $S$ and ${\bf A}^+_c$ are in one-to-one correspondence ({\bf Revuz correspondence}) 
as follows:
for each smooth measure $\mu$, there exists a unique PCAF
$\{A_{t}\}_{t\geq 0}$
such that 
for any $f\in {\mathscr B}^+(X)$ and $\gamma$-excessive
function $h$ ($\gamma\geq 0$),
that is, $e^{-\gamma{t}}p_{t}h\leq h$, 
\begin{equation}\label{eq:7}
 \lim_{t\rightarrow 0}\frac{1}{t}E_{h\cdot{m}}
  \left(\int_{0}^{t}f(X_{s})dA_{s}\right)
  =\int_{X}f(x)h(x)\mu(dx)
\end{equation}
(\cite[Theorem 5.1.7]{FOT}). 
Here, $E_{h\cdot{m}}(\, \cdot\, )=\int_{X}E_{x}(\, \cdot\, )h(x)m(dx)$.
 We denote by $A^{\mu}_t$ the PCAF corresponding to $\mu\in S$.
For a signed smooth measure $\mu=\mu^{+}-\mu^{-}$, we define
$A^{\mu}_t=A^{\mu^{+}}_t-A^{\mu^{-}}_t$.

We introduce some classes of smooth measures.

\begin{defn}\label{def-Kato}\rm
Suppose that $\mu\in S$ is a positive Radon measure.
\begin{enumerate}[(1)]
	\item A measure $\mu$ is said to be in the {\it Kato class} of $\bfM$
	($\cK$ in abbreviation)
	if 
	\begin{equation*}
	\lim_{\alpha \to \infty} \| G_{\alpha} \mu \|_{\infty} = 0.
	\end{equation*}
	A measure $\mu$ is said to be in the {\it local Kato class} ($\cK_{loc}$ in abbreviation)
	if for any compact set $K$,  $1_K\cdot \mu$ belongs to $\cK$.	
	\item Suppose that $\bfM$ is transient. A measure $\mu$ is said
	to be in the class $\cK_{\infty}$ if for any $\eps > 0$, 
	there exists a compact set $K = K(\eps)$ 
	\begin{equation*}
	\sup_{x \in X} \int_{K^c} G(x,y)\mu(dy) < \eps.
	\end{equation*}
	A measure $\mu$ in $\cK_{\infty}$ is called {\it Green-tight}. 
%	for all measurable sets $B \subset K$ with $\mu(B) < \del$. 
%	\item Suppose that $\bfM$ is recurrent. A measure $\mu$ is said
%	to be in the class $\cK_{\infty}^1$ 
%	if $\mu$ satisfies the property in the definition (2) by replacing $G(x,y)$ by $G_1(x,y)$. 
\end{enumerate}
\end{defn}
We note that every measure treated in this paper is supposed to be Radon. 
%Thus we see from \cite[Theorem 3.9]{ABM} that $\mu\in\cK$ if and only if
%\bequ\la{a-b-m}
%\lim_{t\downarrow 0}\sup_{x\in X}\bf\bfE_x(A^\mu_t)=\lim_{t\downarrow 0}\sup_{x\in X}\int_0^t\int_Xp(s,x,y)\mu(dy)ds=0.
%\eequ
We denote the Green-tight class by $\cK_{\infty}(G)$ if we would like to emphasize 
 the dependence of the Green kernel. Chen \cite{Chen-gauge} define the Green-tight class in slightly
 different way; however the two definitions are equivalent under {\textbf (SF)} (\cite[Lemma 4.1]{KK-0}). 

Let $\mu=\mu^+-\mu^-\in\cK_{loc}-\cK$.
 We define the Schr\"odinger form by

\begin{equation}\la{sch-form}
\left\{
\begin{split}
& \cE^{\mu} (u,u) = \cE(u,u) + \int_X u^2 d\mu  \\
& \cD(\cE^{\mu}) = \cD(\cE) \cap L^2(X;\mu^+).
\end{split}
\right.
\end{equation}
 Denoting by $\cL^{\mu} = \cL - \mu$ the self-adjoint 
 operator generated by the closed symmetric form $(\cE^\mu,\cD(\cE^\mu))$, $(-\cL^{\mu}u,v)_m=\cE^\mu(u,v)$,
we see that the associated semigroup $\exp(t\cL^{\mu})$ is expressed as 
$\exp(t\cL^{\mu})f(x)=\bfE_x\left(\exp(-A^{\mu}_t)f(X_t)\right)$
(cf. \cite{ABM}). 

Let ${\bf M}^{\mu^+}=(P^{\mu^+}_x,X_t,\zeta )$ the subprocess of ${\bf M}$ by the multiplicative 
functional $\exp({-A^{\mu^+}_t})$ and suppose that ${\bf M}^{\mu^+}$ is also strong Feller (For this 
conditions, refer \cite{CK}). 

\section{Maximum Principle}
In this section we consider the maximum principle for Schr\"odinger forms.
For $h\in\sB(X)$ we denote by $h^+$ and $h^-$ the positive and negative part of $h$.

\begin{thm}\la{MP}
Assume {\rm ({\bf A})}. Then
$$
\lambda (\mu )>1\Longleftrightarrow {({\bf MP}) }.
$$
\end{thm}
\begin{proof}
For $h\in\cH^{ba}(\mu)$
\begin{align*}
h(x)& \leq \bfE_x\left(e^{-A^{\mu}_t}h(X_t)\right)= \bfE_x^{\mu^+}\left(e^{A^{\mu^-}_t}h(X_t)\right)\leq \bfE_x^{\mu^+}\left(e^{A^{\mu^-}_\zeta}h^+(X_t)\right)\\
& \leq \|h^+\|_\infty \cdot \bfE_x^{\mu^+}\left(e^{A^{\mu^-}_\zeta};t<\zeta \right).
\end{align*}
If $\lambda (\mu )>1$, then $\sup_{x\in X}\bfE^{\mu^+}_x(\exp(A^{\mu^-}_\zeta ))<\infty $
by \cite[Theorem 5.1]{Chen-gauge}. Hence 
the right-hand side tends to 0 as $t\to\infty $ 
because 
$$
\lim_{t\to\infty }\bfP^{\mu^+}_x(t<\zeta )=\bfE_x\left(e^{-A^{\mu^+}_\infty}1_{\{\zeta=\infty\} }\right)=0
$$
by Assumption ({\bf A}). 

Suppose $\lambda (\mu )\leq 1$. By the definition of $\lambda (\mu )$
\bequ\la{op}
 \inf\left\{\cE^{\mu^+}(u,u)\mid\lambda (\mu )\int_Xu^2d\mu^-=1\right\}=1.
\eequ
It follows from \cite[Lemma 5.16, Corollary 5.17]{T3} that the minimizer $h$ in (\ref{op}) is  
a bounded  positive continuous with
$p^{\mu^+-\lambda (\mu)\mu^-}_t$-invariance, $h(x)=p^{\mu^+-\lambda (\mu)\mu^-}_th(x)$. Hence  
$$
h(x)=p^{\mu^+-\lambda (\mu)\mu^-}_th(x)\leq p^{\mu^+-\mu^-}_th(x)=p^{\mu}_th(x),
$$
and ({\bf MP}) does not hold. 
\end{proof}

\medskip
In the sequel of this section, we deal with a strongly local Dirichlet form and extend a result of \cite{BNV}.
 We set 
\begin{align*}
{\cS}&=\left\{\{x_n\}_{n=1}^\infty \subset X\mid
\text{$x_n\to \triangle $ and $\lim_{n\to\infty }\bfE_{x_n}\left(e^{-\zeta}\right)=1$.}
\right\},\\
\widetilde{\cS}&=\left\{\{x_n\}_{n=1}^\infty\subset X\mid
\text{$x_n\to \triangle $ and $\lim_{n\to\infty }\bfP_{x_n}(\zeta >\epsilon )\to 0$ for any $ \epsilon >0$}
\right\}.
\end{align*}

\begin{lem}\la{S=S} It holds that  
$$
\cS=\widetilde{\cS}.
$$
\end{lem}
\begin{proof}
For $\{x_n\}_{n=1}^\infty\in {\cS}$
$$
\bfE_{x_n}(e^{-\zeta })\leq e^{-\epsilon }\bfP_{x_n}(\zeta >\epsilon )+\bfP_{x_n}(\zeta\leq \epsilon )=1-(1-e^{-\epsilon })\bfP_{x_n}(\zeta > \epsilon ),
$$
and thus
$$
\varlimsup_{n\to\infty }\bfP_{x_n}(\zeta > \epsilon )\leq \varlimsup_{n\to\infty }\frac{1-E_{x_n}(e^{-\zeta })}{1-e^{-\epsilon }}=0.
$$

For $\{x_n\}_{n=1}^\infty\in \widetilde{\cS}$
%\begin{align*}
$$
\bfE_{x_n}(e^{-\zeta })=\bfE_{x_n}(e^{-\zeta };\zeta >\epsilon )+\bfE_{x_n}(e^{-\zeta };\zeta\leq \epsilon )
\geq e^{-\epsilon }\bfP_{x_n}(\zeta \leq \epsilon ),
$$
%\end{align*}
and thus 
$\varliminf_{n\to\infty }\bfE_{x_n}(e^{-\zeta })\geq e^{-\epsilon }$ and
$\lim_{n\to\infty }\bfE_{x_n}(e^{-\zeta })=1$.
\end{proof}

A Dirichlet form $(\cE,\cD(\cE))$ is said to be {\it strongly local}, if $\cE(u,v)=0$ for any 
$u,v\in\cD(\cE)$ such that $u$ is constant on a neighborhood of 
$\supp[v]$.
In the sequel of this section, we assume that $(\cE,\cD(\cE))$ is {strongly local}. We introduce
\begin{equation*}
\widetilde{\cH}^{ba}(\mu)=\left\{h\mid
\begin{split}
&\text{$h\in\cD_{loc}(\cE)\cap C(X)$ is bounded above, $\cE^\mu(h,\varphi )\leq 0$ for any}\\
&\text{$\varphi \in \cD(\cE)\cap C^+_0(X)$, $\varlimsup_{n\to\infty }h(x_n)\leq 0$ for any $\{x_n\}_{n=1}^\infty \in \cS$}.
\end{split}
\right\}.
\end{equation*}
%$$
%\cH^{bb}(\mu)= \{ h \in \cD_{loc}(\cE) \cap C(X) \mid h\ \text{is bounded below},\  p_t^{\mu} h \leq  h\}.
%$$

\begin{lem}\la{subseq}
Let $\{\tau_n\}_{n=1}^\infty $ be a sequence of stopping times such that $\tau_n<\zeta $ and $\tau_n\uparrow \zeta $, 
 as $n\to\infty $, $\bfP_x$-a.s. Then there exists a subsequence $\{\sigma _n\}_{n=1}^\infty $ of 
$\{\tau_n\}_{n=1}^\infty $ such that 
\bequ\la{subs}
\bfP_x\left(\{X_{\sigma _n}\}\in \cS\right)=1.
\eequ
\end{lem}
\begin{proof}
First note 
\begin{align*}
\{\zeta (\theta _{\tau _n})>\epsilon ,\ \tau_n<\zeta \}&
=\{\tau_n+\zeta (\theta _{\tau _n})>\tau_n+\epsilon ,\ \tau_n<\zeta \}\\
&=\{\zeta>\tau_n+\epsilon ,\ \tau_n<\zeta \}=\{\zeta>\tau_n+\epsilon \}.
\end{align*}
We then have by the strong Markov property
\begin{align*}
\bfE_x(\bfP_{X_{\tau_n}}(\zeta >\epsilon ))&=\bfE_x(\bfP_{X_{\tau_n}}(\zeta >\epsilon );\tau_n<\zeta )\\
&=\bfE_x(\bfP_x(\zeta (\theta _{\tau_n})>\epsilon,\ \tau_n<\zeta|{\mathscr
F}_{\tau_n}) )\\
&=\bfP_x(\zeta >\tau_n+\epsilon)\longrightarrow 0\ \ \text{as}\  n\to\infty. 
\end{align*}
Hence there exists a subsequence $\{\tau^{(1)}_n\}_{n=1}^\infty $ of $\{\tau_n\}_{n=1}^\infty $
 such that
$$
\bfP_{X_{\tau^{(1)}_n}}(\zeta >1)\longrightarrow 0\ \ \text{as}\ \ n\to\infty,\ \ \bfP_x\text{-a.s.}.
$$
By the same argument 
\begin{align*}
\bfE_x(\bfP_{X_{\tau^{(1)}_n}}(\zeta >1/2);\tau^{(1)}_n<\zeta )\longrightarrow 0
\end{align*}
and there exists a subsequence $\{\tau^{(2)}_n\}_{n=1}^\infty$ of $\{\tau^{(1)}_n\}_{n=1}^\infty$ such that
$$
\bfP_{X_{\tau^{(2)}_n}}(\zeta >1/2)\longrightarrow 0\ \ \text{as}\ \ n\to\infty,\ \ \bfP_x\text{-a.s. }
$$
By continuing this procedure we can take a subsequence $\{\tau^{(k)}_n\}_{n=1}^\infty$ of 
$\{\tau^{(k-1)}_n\}_{n=1}^\infty$ such that
$$
\bfP_{X_{\tau^{(k)}_n}}(\zeta >1/k)\longrightarrow 0\ \ \text{as}\ \ n\to\infty,\ \ \bfP_x\text{-a.s. }
$$
The sequence $\{\sigma _n:=\tau^{(n)}_n\}_{n=1}^\infty $ is a desired one.
\end{proof}

\begin{lem}\la{meas} Suppose $(\cE,\cD(\cE))$ is strongly local. 
Let $\{D_n\}_{n=1}^\infty$ be a sequence of relatively compact open sets such that $D_n\uparrow X$.
Define $S_n=\inf\{t>0\mid A^{\mu^-}_t>n\}$ and $T_n=S_n\wedge \tau_{D_n}$. 
Then for $h\in\widetilde{\cH}^{ba}(\mu)$
$$
\bfE_x\left(e^{-A^{\mu}_{T_n\wedge t}}h(X_{T_n\wedge t})\right)\geq h(x)\ \ \rm{q.e.}\ x.
$$
\end{lem}

\begin{proof}
This lemma can be derived by the argument similar to that in \cite[Lemma 4.7]{T3}. In fact, 
 put ${\mathcal L}=\cD(\cE)\cap C_0(X)$. Then ${\mathcal L}$ 
 is a {\it Stone vector lattice}, i.e., if $f,g\in {\mathcal L}$,
 then $f\vee g\in {\mathcal L}$, $f \wedge 1\in {\mathcal L}$.
 For $h\in \widetilde{\cH}^{ba}(\mu)$ define the functional $I$  
 by
 \begin{equation}\la{i-r}
  I(\varphi)=-\cE^{\mu}(h,\varphi ), \quad \varphi \in{\mathcal L}. 
 \end{equation}
Then $I(\varphi )$ is a pre-integral, that is, 
$I(\varphi_{n})\downarrow 0$ whenever $\varphi_{n}\in{\mathcal L}$ and $\varphi_n(x)\downarrow 0$
 for all $x\in X$. Indeed, let $\psi \in\cD(\cE)\cap C_0(X)$ such that $\psi =1$ on supp[$\varphi _1$].
 Then $\varphi _n\leq \|\varphi_n\|_\infty \psi$ and 
$$
I(\varphi _n)\leq \|\varphi_n\|_\infty\cdot  I(\psi )\downarrow 0, \ \ n\rightarrow \infty .
$$  
Notice that by the regularity of $(\cE,\cD(\cE))$ 
the smallest $\sigma $-field generated by ${\mathcal L}$ is identical with 
the Borel $\sigma $-field. We then see from \cite[Theorem 4.5.2]{Dud} that there exists a positive Borel measure
$\nu$ such that  
\begin{equation}\la{i-def}
I(\varphi )=\int_X\varphi d\nu.
\end{equation}

By the definition of $\nu$ we see that 
$\nu$ is a Radon measure and satisfies (S2) for any increasing sequence $\{F_n\}$ of compact sets with $F_n\uparrow X$.
Let $K$ be a compact set of zero capacity. Then for a
relatively compact open set $D$ such that $K\subset D$, there
exists a sequence $\{\varphi _n\} \subset \cD(\cE)\cap
C^+_0(D) $ such that $\varphi _n\geq 1$ on $K$ and  
$\cE_1(\varphi _n,\varphi _n)\to 0$ as $n\to\infty $ (\cite[Lemma 2.2.7]{FOT}).
For $\psi \in \cD(\cE)\cap C_0(X)$ with $\psi =1$ on $D$,
\begin{equation*}
I(\varphi _n)=-\cE^\mu(h,\varphi _n)=-\cE^\mu(h\psi ,\varphi _n)\leq
 \cE^{|\mu|}(h\psi ,h\psi)^{1/2}\cdot 
 \cE^{|\mu|}(\varphi _n,\varphi _n)^{1/2},
 \end{equation*}
 where $|\mu|=\mu^++\mu^-$. Note that $1_D|\mu|\in\cK$ and $\|G_1(1_D|\mu|)\|_\infty<\infty $. 
 We then see from the Stollmann-Voigt inequality (\cite{Stol-Vo}) that
\begin{equation*}
\int_X\varphi _n^2d|\mu|=\int_X\varphi _n^21_Dd|\mu|\leq \|G_1(1_D|\mu|)\|_\infty \cdot \cE_1(\varphi
 _n,\varphi _n)\longrightarrow 0, \quad n\to\infty 
\end{equation*}
and $\cE^{|\mu|}(\varphi _n,\varphi _n)\rightarrow 0$ as $n\to \infty $.
Since  
\begin{equation*}
\nu(K)\leq \int_X\varphi _nd\nu= I(\varphi _n)\to 0,\ \  n\to\infty
\end{equation*}
$\nu$ satisfies (S1), consequently the measure $\nu$ is smooth. 

The equations (\ref{i-r}), (\ref{i-def}) lead us to 
\begin{equation*}
\cE(h,\varphi )=-\int_X\varphi hd\mu-\int_X\varphi d\nu=-\int_X\varphi (hd\mu+d\nu).
\end{equation*}
On account of \cite[Theorem 5.5.5]{FOT}, we have 
\begin{equation*}
h(X_t)=h(X_0)+M^{[h]}_t+\int_0^th(X_s)dA^{\mu}_s+A^{\nu}_t \ \ \ \bfP_x\text{-a.s., q.e. x}.
\end{equation*}
Hence, by It$\hat{\rm o}$'s formula 
\begin{equation*}
 \begin{split}
  e^{-A^\mu_t}h(X_t)&=
  h(X_0)+\int_{0}^{t}e^{-A^{\mu}_{s}}dM_s^{[h]}
  +\int_{0}^{t}e^{-A^{\mu}_{s}}h(X_{s})dA^{\mu}_{s}\\
&\ \ +\int_0^te^{-A^\mu_s}dA^\nu_s-\int_0^te^{-A^\mu_s}h(X_s)dA^\mu_s\\
  &=h(X_0)+\int_0^te^{-A^\mu_s}dM_s^{[h]}+\int_0^te^{-A^\mu_s}dA^\nu_s \ \ \ \bfP_x\text{-a.s., q.e. x}.
 \end{split}
\end{equation*}
Since $\int_0^{T_n\wedge t}e^{-A^\mu_s}dM_s^{[h]}$ is a martingale and $\int_0^te^{-A^\mu_s}dA^\nu_s\geq 0$,
 \begin{equation*}
  E_{x}\left(e^{-A^{\mu}_{T_n\wedge t}}h(X_{T_n\wedge t})\right)\geq h(x) \ \ \text{q.e. x}.
 \end{equation*}
\end{proof}

\begin{lem}\la{include} It holds that
$$
\widetilde{\cH}^{ba}(\mu)\subset {\cH}^{ba}(\mu).
$$
\end{lem}
\begin{proof}
Let $h$ be a function in $\widetilde{\cH}^{ba}(\mu)$ and $\{T_n\}_{n=1}^\infty$ 
a sequence of stopping times defined in Lemma \ref{meas}. 
We fix a point $x\in X$ such that 
$$
\bfE_x\left(e^{-A^{\mu}_{T_n\wedge t}}h(X_{T_n\wedge t})\right)\geq h(x).
$$
Since $T_n<\zeta $ and
$T_n\uparrow \zeta $, 
%By Lemma \ref{meas}, we have
%$$
%\bfE_x\left(e^{-A^{\mu}_{\sigma _n\wedge t}}h(X_{\sigma_n\wedge t})\right)\geq h(x)\ \ \text{q.e. x}.
%$$
 we can take a subsequence $\{\sigma_n\}$ of $\{T_n\}$ satisfying (\ref{subs}) 
in Lemma \ref{subseq}. Since $h^+$ is bounded continuous and 
$$
\varlimsup_{n\to\infty }h^+(X_{\sigma_n\wedge t})=0,\ \ \bfP_x\text{-a.s. on}\ \{t\geq \zeta\} 
$$
by (\ref{subs}), we have
\begin{align*}
&\ \ \ \ \varlimsup_{n\to\infty }\bfE_x\left(e^{-A^{\mu}_{\sigma_n\wedge t}}h^+(X_{\sigma_n\wedge t})\right)\\
&\leq \varlimsup_{n\to\infty }\bfE_x\left(e^{-A^{\mu}_{\sigma_n\wedge t}}h^+(X_{\sigma_n\wedge t});t<\zeta \right)+\varlimsup_{n\to\infty }\bfE_x\left(e^{-A^{\mu}_{\sigma_n\wedge t}}h^+(X_{\sigma_n\wedge t});t\geq \zeta \right)\\
&\leq \bfE_x(\left(e^{-A^{\mu}_{t}}h^+(X_{t});t<\zeta \right)+\bfE_x\left(\varlimsup_{n\to\infty }e^{-A^{\mu}_{\sigma_n\wedge t}}h^+(X_{\sigma_n\wedge t});t\geq \zeta \right)\\
&= \bfE_x\left(e^{-A^{\mu}_{t}}h^+(X_{t}) \right).
\end{align*}
Here, the second inequality above follows from the inverse Fatou's lemma because
$$
e^{-A^{\mu}_{\sigma_n\wedge t}}h^+(X_{\sigma_n\wedge t})\leq e^{A^{\mu^-}_{\sigma_n\wedge t}}h^+(X_{\sigma_n\wedge t})\leq \|h^+\|_\infty\cdot e^{A^{\mu^-}_t}\in
L^1(\bfP_x)
$$
by $\mu^-\in\cK$.

Besides, we have
\begin{align*}
&\ \ \ \ \varliminf_{n\to\infty }\bfE_x\left(e^{-A^{\mu}_{\sigma_n\wedge t}}h^-(X_{\sigma_n\wedge t})\right)\\
&\geq  \varliminf_{n\to\infty }\bfE_x\left(e^{-A^{\mu}_{\sigma_n\wedge t}}h^-(X_{\sigma_n\wedge t});t<\zeta \right)+\varliminf_{n\to\infty }\bfE_x\left(e^{-A^{\mu}_{\sigma_n\wedge t}}h^-(X_{\sigma_n\wedge t});t\geq \zeta \right)\\
&\geq \bfE_x\left(\varliminf_{n\to\infty }e^{-A^{\mu}_{\sigma_n\wedge t}}h^-(X_{\sigma_n\wedge t});t< \zeta \right)= \bfE_x\left(e^{-A^{\mu}_{t}}h^-(X_{t})\right).
\end{align*}
Hence
\begin{align*}
h(x)&\leq \varliminf_{n\to\infty }\bfE_x\left(e^{-A^{\mu}_{\sigma_n\wedge t}}h(X_{\sigma_n\wedge t})\right)\\
&\leq \varlimsup_{n\to\infty }\bfE_x\left(e^{-A^{\mu}_{\sigma_n\wedge t}}h^+(X_{\sigma_n\wedge t})\right)-\varliminf_{n\to\infty }\bfE_x\left(e^{-A^{\mu}_{\sigma_n\wedge t}}h^-(X_{\sigma_n\wedge t})\right)\\
&\leq \bfE_x\left(e^{-A^{\mu}_{ t}}h^+(X_{t}) \right)-\bfE_x\left(e^{-A^{\mu}_{ t}}h^-(X_{t}) \right)= \bfE_x\left(e^{-A^{\mu}_{ t}}h(X_{t})\right),
\end{align*}
and $h(x)\leq p^\mu_th(x)$ for q.e. $x$. 
Since $p^{\mu}_t$ is strong Feller, $p^{\mu}_t(h\vee (-n))(x)\geq h(x)$ for all $x\in X$ and 
$p^{\mu}_th(x)\geq h(x)$ for all $x\in X$ by letting $n$ to $\infty $.
\end{proof}

Following \cite{BNV}, we define the {\it refined maximum principle}:

\bigskip
\noindent
\ \ \ \ (${\bf{RMP}}$)\ \ If $h\in \widetilde{\cH}^{ba}(\mu)$, then $h(x)\leq 0$ for all $x\in X$.

\bigskip
Combining Lemma \ref{include} with Theorem \ref{MP}, we have the next theorem.

\begin{thm}\la{RMP} Suppose $(\cE,\cD(\cE))$ is strongly local. Then under Assumption {\rm{({\bf A})}}
$$
\lambda (\mu)>1 \Longrightarrow ({\bf RMP}).
$$
\end{thm}

\bigskip
\begin{rem}\la{S=S1}\rm
Suppose $D$ is a bounded domain in ${\mathbb R}^d$ and consider the absorbing Brownian motion $(\bfP_x,B_t,\tau_D)$ on $D$, 
where $\tau_D$ is the first exit time from $D$.
If $D$ is Green-bounded, i.e., $\sup_{x\in D}\bfE_x(\tau_D)<\infty $, then $\cS$ is identical to
the set of sequences $\{x_n\}$ such that $x_n\to\partial D$ and 
$E_{x_n}(\tau_D)\to 0$ as $n\to\infty $. Indeed, take $\delta >0$ so that $\sup_{x\in D}\bfE_x(\delta \tau_D)<1$. 
Then since $\sup_{x\in D}\bfE_x(\exp(\delta \tau_D))<\infty $
by Has'minskii's lemma, we see $\sup_{x\in D}\bfE_x(\tau_D^2)<\infty $. 
Hence if $\bfP_{x_n}(\tau_D>\epsilon )\to 0$ as $n\to\infty $, then  
$$
E_{x_n}(\tau_D)\leq  \left(E_{x_n}(\tau_D^2)\right)^{1/2}\cdot  \left(\bfP_{x_n}(\tau_D>\epsilon )\right)^{1/2}+\epsilon 
\bfP_{x_n}(\tau_D\leq \epsilon )\to \epsilon .
$$
Since $\epsilon \bfP_{x_n}(\tau_D> \epsilon )\leq E_{x_n}(\tau_D)$,
the converse follows from Lemma \ref{S=S}.

Let 
$$
\lambda _0=\inf\left\{\frac{1}{2}{\mathbb D}(v,v)+\int_Dv^2d\mu\mid v\in H_0^1(D), \int_Dv^2dx=1\right\},
$$
where ${\mathbb D}$ is the classical Dirichlet integral. 
We see from \cite[Theorem 1.1]{BNV} that
$$
\lambda _0>0\Longleftrightarrow ({\bf RMP}).
$$
Moreover, we see from Lemma \ref{first} below that if $D$ is bounded, then $\lambda _0>0$ and $\lambda (\mu )>1$ 
are equivalent, and so
$$
\lambda (\mu )>1\Longleftrightarrow ({\bf RMP}).
$$

We remark that $\lambda _0>0$ implies $\lambda (\mu )>1$ for a general domain $D$ (Lemma \ref{first} below), while $\lambda(\mu)>1$ does not 
imply $\lambda_0>0$ in general. 
In fact, consider $\cL u=(1/2)u''-\mu u$ ($\mu=\alpha \delta _{-1}-\beta \delta _{1},\ \alpha >0,\ 
\beta >0$) on ${\mathbb R}^1$. We define
$$
\lambda ({\alpha ,\beta }):=\lambda (\mu)=\inf\left\{\frac{1}{2}{\mathbb D}(u,u)+\alpha u({-1})^2\mid u\in H^1(\bR^1),\ \beta u(1)^2=1
\right\}
$$
and
$$
\lambda_0 ({\alpha ,\beta }) :=\inf\left\{\frac{1}{2}{\mathbb D}(u,u)+\alpha u({-1})^2-\beta u(1)^2\mid u\in H^1(\bR^1),\ \int_{{\mathbb R}^1}u^2dx=1
\right\}.
$$
Denote by $\cL_0$ the operator $1/2(d^2/dx^2)-\alpha\delta_{-1}$.
By the Dirichlet principle, the infimum of $\lambda ({\alpha ,\beta })$ is attained by 
the $\cL_0$-harmonic function $u_0$ with $u_0(1)=1/\sqrt{\beta}$, i.e.,
\begin{equation}\la{sch}
u_0(x)=\left\{
\begin{split}
& \gamma,  & x\leq -1,  \\
& \gamma+\frac{{1/\sqrt{\beta}-\gamma}}{2}(x+1), & -1\leq x<1, \\
& 1/\sqrt{\beta}, & x\geq 1.\nonumber
\end{split}
\right.
\end{equation}
Here, $\gamma$ is determined by 
$$
\cL_0u_0(-1)=0 \Longleftrightarrow \frac{u_0'(-1+)-u_0'(-1-)}{2}=\alpha u_0(-1)\Longleftrightarrow \frac{{1/\sqrt{\beta}-\gamma}}{4}=\alpha\gamma,
$$
and thus $\gamma={1}/({\sqrt{\beta}(4\alpha+1)})$. Note that $u_0$ belongs to the extended Dirichlet space $H_e^1(\bR^1)(\supset H^1(\bR^1))$ (cf.
\cite[Exercise 6.4.9]{FOT}). We then see that
$$
\lambda ({\alpha ,\beta })=\frac{1}{2}\int_{-1}^1\left(\frac{du_0}{dx}\right)^2dx+\alpha u_0(-1)^2
=\frac{\alpha}{\beta(4\alpha+1)}.
$$

For $\beta <1/4$, let $\alpha _0={\beta }/{(1-4\beta) }$. Then   
$\lambda (\alpha_0 ,\beta )=1$  and 
$\lambda ({\alpha ,\beta })>1$ for $\alpha > \alpha_0$.
We see from \cite[Lemma 2.2]{TT} that 
$\lambda (\alpha ,\beta )\geq 1$ is equivalent with $\lambda_0 (\alpha,\beta )\geq 0$.  
Noting that $\lambda_0 ({\alpha ,\beta })\leq 0$ for any $\alpha ,\ \beta $, we see that 
 for $\beta <1/4$ and $\alpha > {\beta }/{(1-4\beta) }$, $
 \lambda_0 ({\alpha ,\beta })=0\ \ \text{and}\ \ \lambda ({\alpha ,\beta })>1.
$
\end{rem}

\begin{lem}\la{first} It holds that 
$$
\lambda _0:=\inf\left\{\cE^{\mu}(u,u)\mid u\in\cD(\cE),\ \int_Xu^2dm=1\right\}>0\Longrightarrow \lambda (\mu )>1.
$$
If there exists a positive constant $C$ such that
$$
\int_Xu^2dm\leq C\cE^{\mu^+}(u,u),
$$
then the converse also holds.
\end{lem}

\begin{proof}
Let $\varphi _0\in \cD_e(\cE^{\mu^+})$ be the minimizer in (\ref{i-ground}):
$$
\lambda (\mu )=\cE^{\mu^+}(\varphi _0,\varphi _0),\ \ 
\int_X\varphi _0^2d\mu^-=1.
$$
If $\lambda _0>0$, then 
$$
\lambda (\mu )-1=\cE^{\mu^+}(\varphi _0,\varphi _0)-\int_X\varphi _0^2d\mu^-=\cE^{\mu}(\varphi _0,\varphi _0)
= \lambda _0\int_X\varphi _0^2dm>0.
$$

If $\lambda (\mu )>1$, then  for any $u\in\cD(\cE)$
$$
\cE^{\mu^+}(u,u)-\lambda(\mu)\int_Xu^2d\mu^-\geq 0\Longleftrightarrow \lambda(\mu)\cdot \cE^{\mu}(u,u)
\geq (\lambda(\mu)-1)\cdot  \cE^{\mu^+}(u,u).
$$
Hence by the assumption,
$$
\cE^{\mu}(u,u)\geq \frac{(\lambda(\mu)-1)}{C\lambda(\mu)} \int_Xu^2dm.
$$

\end{proof}

\section{Liouville Property}
Let us introduce the set of $p^\mu_t$-invariant bounded functions by
$$
\ \cH^b(\mu)= \{ h\in{\mathscr B}_b(X)  \mid \  p_t^{\mu} h = h\}.
$$
We here define the Liouville property ({\bf L}) by

\bigskip
\noindent
\ \ \ \ (${\bf L}$)\ \ If $h\in {\cH}^{b}(\mu)$, then $h(x)=0$ for all $x\in X$.

\bigskip
\begin{cor}\la{Cor} {Suppose {\rm({\bf A})}. Then}
$$
\lambda (\mu )>1\Longrightarrow ({\bf L}).
$$
\end{cor}
\begin{proof}
Let
$$
\cH^{bb}(\mu)= \{ h\in\sB(X) \mid h\ \text{ is bounded below},\  p_t^{\mu} h \leq  h\}.
$$
We see, by the same argument as in Theorem \ref{MP}, that an element $h$ in $\cH^{bb}(\mu)$ satisfies
$h(x)\geq 0$ for any $x\in X$. 
Since $ \cH^b(\mu)= \cH^{ba}(\mu)\cap \cH^{bb}(\mu)$, this corollary is derived.
\end{proof}

For a strongly local Dirichlet form $(\cE,\cD(\cE))$ we set
\begin{equation*}
\widetilde{\cH}^b(\mu)=\{h\in\cD_{loc}(\cE)\cap C_b(X)\mid
\cE^\mu(h,\varphi )=0, \forall \varphi \in \cD(\cE)\cap C_0(X)\}.
\end{equation*}

\begin{lem}\la{con}
Assume $(\cE,\cD(\cE))$ is strongly local. If ${\bfM}$ is conservative, then $\widetilde{\cH}^{b}(\mu)\subset {\cH}^{b}(\mu)$.
\end{lem}
\begin{proof}
For $h\in \widetilde{\cH}^{b}(\mu)$ let $\{T_n\}_{n=1}^\infty $ be the sequence of stopping times 
defined in Lemma \ref{meas}. Then $\bfE_x(\exp(-A^{\mu}_{T_n\wedge t})h(X_{T_n\wedge t}))=h(x)$ for any $n$.
Noticing that $T_n\to \infty $, $\bfP_x$-a.s. 
by the conservativeness of ${\bf M}$ and that $\exp({-A^{\mu}_{T_n\wedge t}})h(X_{T_n\wedge t})
\leq \|h\|_\infty\exp(A^{\mu^-}_t)\in L^1(\bfP_x)$, we have 
$$
h(x)=\lim_{n\to\infty }\bfE_x\left(e^{-A^{\mu}_{T_n\wedge t}}h(X_{T_n\wedge t})\right)
=\bfE_x\left(e^{-A^{\mu}_t}h(X_t)\right)
$$
by the dominated convergence theorem.
\end{proof}

\bigskip
Define the property $(\widetilde{\bf L})$ by

\bigskip
\noindent
\ \ \ \ $(\widetilde{\bf L})$\ \ If $h\in \widetilde{\cH}^{b}(\mu)$, then $h(x)=0$ for all $x\in X$.

\bigskip
\noindent
Lemma \ref{con} leads us to the next corollary.

\begin{cor}\la{Cor1} Suppose $(\cE,\cD(\cE))$ is strongly local and ${\bfM}$ is conservative. 
Then under Assumption {\rm({\bf A})}
$$
\lambda (\mu )>1\Longrightarrow (\widetilde{\bf L}).
$$
\end{cor}

\bigskip
%Consider the Brownian motion $(\bfP_x,B_t)$ on $\bR^d$, $d\geq 3$.
%Pinsky \cite{RP} treat absolutely continuous potentials $d\mu=V^+dx-V^-dx$ and prove in  \cite[Theorem 1.1]{RP}
%that if $
%\sup_{x\in{\mathbb R}^d}\bfE_x\left(\exp({\int_0^\infty V^-(B_t)dt})\right)<\infty,
%$
%the Liouville property $(\widetilde{\bf L})$ is equivalent to 
%$$
%\int_0^\infty V^+(B_t)dt=\infty,\ \bfP_x\text{-a.e.}\ (\Longleftrightarrow \ ({\bf A})).
%$$ 

We finally give a Schr\"odinger operator, $-1/2\Delta +\mu$ which satisfies  $(\widetilde{\bf L})$; however, the positive part and negative part of potential $\mu$ satisfy
$$
\bfE_x\left(e^{-A^{\mu^+}_\infty}\right)=0,\ \ \ \sup_{x\in{\mathbb R}^d}\bfE_x\left(e^{A^{\mu^-}_\infty}\right)=\infty.
$$ 

\begin{example} \la{e-2}\rm
Let us define 
$$
\lambda _1=\inf\left\{\frac{1}{2}{\mathbb D}(u,u)\mid u\in H^1({\mathbb R}^d),\int_{{\mathbb R}^d}u^2d\sigma =1\right\}
$$
and
$$
\lambda _2=\inf\left\{\frac{1}{2}{\mathbb D}(u,u)+(u,u)_m\mid u\in H^1({\mathbb R}^d),\int_{{\mathbb R}^d}u^2d\sigma =1\right\},
$$
where $m$ is the Lebesgue measure and $\sigma $ the measure such that $\sigma|_{\partial B(0,1)}$
 is the surface measure of $\partial B(0,1)$ and $\sigma(\bR^d\setminus \partial B(0,1))=0$. 
Let $\mu=m-\gamma \sigma $, that is, $\mu^+=m$,  $\mu^-=\gamma \sigma $ ($\gamma>0$). 
Note that $A^{m}_t=t$ and $A^\sigma _t$ is the local time of the unit sphere. 
We see that if $\lambda_1< \gamma  <\lambda_2$, then $\lambda (\mu)>1$, and $-1/2\Delta +\mu$ satisfies ($\widetilde{\bf L}$); however,
$
\bfE_x\left(\exp\left(A^{\gamma \sigma} _\infty \right)\right)=\infty .
$
\end{example} 
\bibliographystyle{amsplain}

\begin{thebibliography}{10}

\bibitem{ABM}
Albeverio, S, Blanchard, P. and  Ma, Z.-M.: 
\newblock Feynman-{K}ac semigroups in terms of signed smooth measures,
\newblock In {\em Random partial differential equations ({O}berwolfach, 1989)},
     Birkh\"auser,
   (1991), 1-31.
 \bibitem{BNV} Berestycki, H., Nirenberg, L., Varadhan, S. R. S.: The principal eigenvalue and maximum principle for second-order elliptic operators in general domains, {\it Comm. Pure Appl. Math.}, {\bf 47}, (1994), 47-92.
 \bibitem{Chen-gauge}Chen, Z.-Q.:
	 Gaugeability and conditional gaugeability, 
	 {\it Trans. Amer. Math. Soc.}, {\bf 354}, (2002), 4639-4679.
 \bibitem{CK} Chen, Z.-Q. and Kuwae, K.: On doubly Feller property, {\it Osaka
	 J. Math.} {\bf 46},  (2009), 909-930. 

 \bibitem{Dud}	Dudley, R. M.: Real Analysis and Probability, Cambridge
	 Studies in Advanced Mathematics, 74, Cambridge University
	 Press, (2002). 
 \bibitem{FOT}Fukushima, M., Oshima, Y. and Takeda, M.:
	 Dirichlet Forms and Symmetric Markov Processes, 
	 Walter de Gruyter, 2nd ed. (2011).
% \bibitem{Ge} Getoor, R. K.: Transience and recurrence of Markov
%	 processes, Lecture Notes in Math., {\bf 784}, Springer, (1980),
%	 397-409.
\bibitem{GH}Grigor'yan, A., Hansen, W.: A Liouville property for Schr\"odinger operators, 
{\it Math. Ann.}, {\bf 312}, (1998), 659-716.
\bibitem{KK-0}Kim, D. and  Kuwae, K.: Analytic characterizations of gaugeability for 
generalized Feynman-Kac functionals, to appear in {\it Trans. Amer. Math. Soc.}
% \bibitem{KK} Kim, D. and  Kuwae, K.: On a stability of heat kernel
%	 estimates under generalized 
%	 non-local Feynman-Kac perturbations for stable-like processes,
%	 preprint. 


 	\bibitem{RP} Pinsky, R.G.: A probabilistic approach to bounded/positive solutions for Schr\"odinger 
operators with certain classes of potentials, {\it Trans. Amer. Math. Soc.}, {\bf 360}, (2008), 6545-6554. 
%\bibitem{schmu} Schmuland, B.: Extended Dirichlet spaces, {\it CRMR, Royal Soc. Canada}, {\bf 21}, (1999), 
%146-152. 
% \bibitem{S}Silverstein, M. L.: Symmetric Markov Processes, Lecture Notes in
% Mathematics {\bf 426}, Springer-Verlag, (1974).
 \bibitem{Stol-Vo}Stollmann, P. and Voigt, J.:
	 Perturbation of Dirichlet forms by measures, {\it Potential
	 Anal.}, {\bf 5}, (1996), 109-138. 
% \bibitem{T1}Takeda, M.: Gaugeability for Feynman-Kac functionals with applications to
%	 symmetric $\alpha$-stable processes, 
%	 {\it Proc. Amer. Math. Soc.}, {\bf 134}, (2006), 2729-2738.
 \bibitem{T2}Takeda, M.: A variational formula for Dirichlet
	 forms and existence of ground states, {\it J. Funct. Anal.},
	 {\bf 266}, (2014), 660-675. 
	  \bibitem{T3}Takeda, M.: Criticality and subcriticality of generalized Schr\"{o}dinger forms, 
 {\it Illinois J.  Math.}, {\bf 58}, (2014), 251-277. 
 
 \bibitem{TT} Takeda, M., Tsuchida, K.: Differentiability of spectral 
functions for symmetric $\alpha $-stable processes, Trans. Amer. Math. Soc. {\bf 359}, 4031-4054 (2007).
\end{thebibliography}

\end{document}